\documentclass[10pt,a4paper]{amsart}
\usepackage[latin1]{inputenc}
\usepackage{amsmath}
\usepackage{amsfonts}
\usepackage{amssymb}
\usepackage{bm}
\usepackage{url} 
\usepackage{tabularx}
\usepackage{proof}
\usepackage{graphicx}
\usepackage{subfig}

\newtheorem{thm}{Theorem}

\newtheorem{prop}[thm]{Proposition}

\title{Approximate Implicitization of Triangular B\'ezier Surfaces}

\author{Oliver J. D. Barrowclough}
\author{Tor Dokken}
\email{oliver.barrowclough@sintef.no, tor.dokken@sintef.no}
\address{SINTEF ICT, Applied Mathematics\\P.O. Box 124, Blindern \\ 0314 Oslo, Norway}

\begin{document}

\begin{abstract}
We discuss how Dokken's methods of approximate implicitization can be applied to triangular B\'ezier surfaces in both the original and weak forms. The matrices $\mathbf{D}$ and $\mathbf{M}$ that are fundamental to the respective forms of approximate implicitization are shown to be constructed essentially by repeated multiplication of polynomials and by matrix multiplication. A numerical approach to weak approximate implicitization is also considered and we show that symmetries within this algorithm can be exploited to reduce the computation time of $\mathbf{M}.$ Explicit examples are presented to compare the methods and to demonstrate properties of the approximations.
\end{abstract}

\maketitle

\section{Introduction}

Methods for conversion between the two main representations of curves and surfaces in CAGD, namely the parametric and implicit forms, have been widely investigated within the CAGD community. Of these, the parametric form has established itself as the representation of choice in most CAGD systems due to its intuitive geometric nature \cite{Hoffman_implicit}. However, the implicit form has distinct advantages over the parametric form in solving certain geometrical problems and thus the possibility to have a dual representation is, in some circumstances, appealing \cite{juttler_approxmethods}. For example, the implicit representation allows us to immediately determine whether a given point lies on the curve or surface. Although exact formulas can be devised for low degree surfaces, higher order parametric geometries require computationally expensive algorithms such as recursive subdivision. Implicit representations are also useful in intersection problems. Notably, ray tracing of implicitly defined surfaces is much quicker than ray tracing of parametric surfaces. Despite these advantages, exact implicit representations of rational parametric curves and surfaces lead to high polynomial degrees in the implicit equation. In general, implicit equations of high degree are not desirable due to often having extraneous branches and singularities that are not necessarily present in the parametric form. They also exhibit a lack of numerical stability \cite{Sederberg_approx}. 

The procedure of converting from the rational parametric to the implicit form of a curve or surface is called implicitization. Traditional methods using Gr\"{o}bner bases or resultants, focused solely on exact implicitization. Exact implicit representations use exact arithmetic, whereas in CAD and CAGD, the use of floating point arithmetic is desirable due to performance. Approximate implicitization provides numerically stable methods to approximate a parametric surface using lower degree implicit equations. In \cite{dokken_thesis}, a method for approximate implicitization was introduced that allows us to choose the degree of the implicit equation to be defined. The theory behind this approach to approximate implicitization of rational parametric manifolds in $\mathbb{R}^l$ has been thoroughly developed in \cite{dokken_thesis,dokken_approx}. A similar approach known as weak approximate implicitization was developed in \cite{dokken_weak}. The aim of this paper is to present both these methods, in the special case of approximate implicitization of triangular B\'ezier surfaces. Although other methods of approximate implicitization exist \cite{pratt_direct,taubin_estimation,Sederberg_approx,wurm_approx}, the original and weak methods that we follow provide fast algorithms with a high order of convergence that are well suited to curves and surfaces defined in a partition of unity basis \cite{thomassen_selfint}. 

While approximate implicitization of tensor-product B\'ezier surfaces is a fairly simple extension of approximate implicitization of 2D rational parametric curves, triangular B\'ezier surfaces are somewhat more difficult. They are, however, expressed naturally in terms of a bivariate Bernstein basis over a triangular domain, which forms a partition of unity. This allows us to follow the steps of original approach fairly directly.


Unfortunately, implicitization algorithms tend to be computationally expensive and as such are hindered in CAGD applications. This paper will highlight some symmetries in the numerical approach to the algorithm that can be exploited to reduce the computation time, and thus improve the prospect of dual representations in CAGD. 


This paper will be organised as follows. Section 2, will briefly introduce the concepts required to define triangular B\'ezier surfaces, and highlight some properties of Bernstein polynomials that are important for approximate implicitization. Section 3 will present the procedure for approximate implicitization in the context of B\'ezier triangles, both in the original and weak forms. It will highlight some new observations that significantly reduce the number of computations required in the numerical form of the algorithm. The accuracy and convergence rates of approximate implicitization will also be stated. Section 4 will describe a simple example of approximate implicitization of a B\'ezier triangle before concluding with some examples that are more relevant in practice.

\section{Triangular B\'ezier Surfaces}

Triangular B\'ezier surfaces, also known as B\'ezier triangles, were developed by Paul de Casteljau to offer a natural generalization of B\'ezier curves to surfaces \cite{farin_triangularbb}. Although, tensor-product patches may be more intuitive (and are certainly used more widely in CAGD), the triangular patches are in some sense a more fundamental generalization. In this section we recall the notation of B\'ezier triangles and state some simple results about Bernstein polynomials. For a comprehensive review of these concepts we refer the reader to \cite{farin_triangularbb,farin_cagd}. 

\subsection{Barycentric Coordinates}

In this paper we will make extensive use of barycentric coordinates, both over triangles and tetrahedra. Barycentric coordinates over triangles provide a natural domain in which to define the B\'ezier triangle, whereas tetrahedral barycentric coordinates will be used to define the implicit surface. We introduce the notation in the general form to capture both these circumstances in a common definition. 

Barycentric coordinates allow us to express any point $\mathbf{x} \in \mathbb{R}^l$ as
\[
\mathbf{x} = \sum_{i=1}^{l+1} \beta_i \mathbf{a}_i, \quad \sum_{i=1}^{l+1} \beta_i = 1,
\]
where $\mathbf{a}_i \in \mathbb{R}^l$ are points defining the vertices of a non-degenerate simplex in $\mathbb{R}^l.$ 

The conversion between Cartesian coordinates $\mathbf{x}=(x_1,\ldots,x_l)$ and barycentric coordinates $\bm{\beta}=(\beta_1,\ldots,\beta_{l+1})$ over the simplex with vertices $(\mathbf{a}_1,\ldots,\mathbf{a}_{l+1}),$ is given by the following relation:
\begin{equation}\label{eq:b2c}
\left( \begin{matrix} \mathbf{x} \\ 1 \end{matrix} \right) = 
\left( \begin{matrix} \mathbf{a}_{1} & \ldots & \mathbf{a}_{l+1} \\
1 & \ldots & 1
 \end{matrix} \right)
 \bm{\beta}.
\end{equation}
If a point lies within the simplex which defines the barycentric coordinate system, the barycentric coordinates of that point are guaranteed to be non-negative. This leads to good numerical stability if all the points in the algorithm are contained within the relevant simplex. We define the domain $\Omega$ to be the triangle formed by a bivariate barycentric coordinate system, and $\Lambda$ to be the tetrahedron formed by a trivariate barycentric coordinate system. Unless explicitly stated, all subsequent coordinates in this paper are assumed to be barycentric.


\subsection{Bernstein Polynomials}

The notation used when describing Bernstein polynomials and B\'ezier triangles is greatly simplified by making use of multi-indices. These provide a natural way to label the basis functions and can be related to regular indices by choosing an ordering. For multi-indices $\mathbf{i} = (i_1,\ldots,i_{l+1})$ and $\mathbf{j} = (j_1,\ldots ,j_{l+1})$ we have the following definitions:
\begin{itemize}
\item $|\mathbf{i}| = i_1 + \cdots + i_{l+1},$
\item $\mathbf{i+j} = (i_1+j_1,\ldots , i_{l+1}+j_{l+1}),$
\item for barycentric coordinates $\bm{\beta},$ define $\bm{\beta}^\mathbf{i}$ $ = \beta_1^{i_1}\cdots\beta_{l+1}^{i_{l+1}},$
\item the multinomial coefficients are defined as 
\begin{displaymath}
\binom{n}{\mathbf{i}} = \frac{n!}{i_1! i_2! \cdots i_{l+1}!},
\end{displaymath}
\item the ordering of choice is the lexicographical ordering, described by $(i_1,\ldots,i_{l+1})<(j_1,\ldots,j_{l+1})$ if and only if there exists an index $k$ such that $i_k<j_k$ and $i_r=j_r$ for all $r<k.$ 
\end{itemize}

We now define the Bernstein basis polynomials of degree $n$ as
\begin{displaymath}
B_{\mathbf{i}}^n(\bm{\beta}) =  \binom{n}{\mathbf{i}} \bm{\beta}^{\mathbf{i}}, \quad |\mathbf{i}| = n,
\end{displaymath}
where $\bm{\beta}$ are barycentric coordinates.

In this paper, care must be taken to distinguish between triangular and tetrahedral Bernstein polynomials as the notation differs only by the variable they are defined under. A triangular Bernstein polynomial will be defined in the variable $\mathbf{s}\in\Omega,$ whereas a tetrahedral Bernstein polynomial will be defined for $\mathbf{u}\in\Lambda.$ We use the variables $\bm{\beta}$ when describing general barycentric coordinates.


We will now state three important properties of Bernstein polynomials that will be used in the implicitization algorithm:

\begin{itemize}
\item The Bernstein basis forms a partition of unity. That is 
\begin{equation}
\sum_{|\mathbf{i}|=n} B_\mathbf{i}^n(\mathbf{\bm{\beta}}) = 1,
\end{equation}
for all barycentric coordinates $\bm{\beta}.$
\item There is a simply derived formula for multiplying Bernstein polynomials of the same form (i.e., triangular or tetrahedral Bernstein polynomials), which is given as follows:
\begin{equation} \label{eq:product}
B_\mathbf{i}^n(\bm{\beta})B_\mathbf{j}^m(\bm{\beta})=\frac{ \binom{n}{\mathbf{i}} \binom{m}{\mathbf{j}} }{ \binom{n+m}{\mathbf{i+j}} }B_\mathbf{i+j}^{n+m}(\bm{\beta}),
\end{equation}
with $|\mathbf{i+j}|=|\mathbf{i}|+|\mathbf{j}|=m+n.$
\item The integral over any Bernstein basis function of given degree is constant. In particular, for the Bernstein basis polynomials over a triangle of unit area \cite{farouki_orthog}: 
\begin{equation}
\int_\Omega{B_\mathbf{i}^n(\mathbf{s})}\ \mathrm{d}\mathbf{s} = \frac{1}{(n+1)(n+2)}. \label{eq:basisint}
\end{equation}
This implies that the integral of any polynomial $q(\mathbf{s})$ defined in the triangular Bernstein basis is given by:
\begin{eqnarray}
\int_\Omega q(\mathbf{s}) \ \mathrm{d}\mathbf{s} = \frac{1}{(n+1)(n+2)} \sum_{|\mathbf{i}|=n} b_{\mathbf{i}}.\label{eq:bernint}
\end{eqnarray}
\end{itemize}

\subsection{B\'ezier Triangles}\label{subsec:bez}

We can now state the definition of a degree $n$ B\'ezier triangle with control points $(\mathbf{c_i})_{|\mathbf{i}|=n}$ in terms of the triangular Bernstein basis as follows:
\begin{equation}
\mathbf{p}(\mathbf{s})=\sum_{|\mathbf{i}|=n} \mathbf{c}_{\mathbf{i}} B_{\mathbf{i}}^n(\mathbf{s}). 
\end{equation}
The control points $\mathbf{c_i}$ can be defined in any space but we restrict them to lie in $\mathbb{R}^3$ since we are interested in surfaces. We consider only the points $\mathbf{s}$ in the domain $\Omega$ so that the entire B\'ezier triangle lies within the convex hull of its control points.

\begin{figure}
\begin{center}
 \begin{tabular}{cc}
 \includegraphics[scale=0.22]{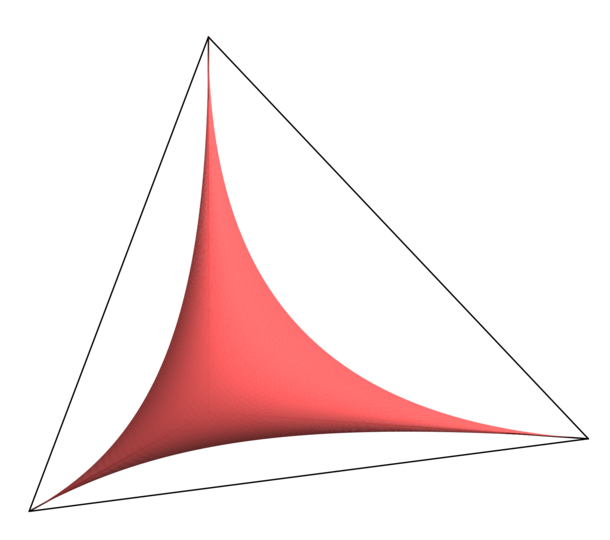} & 
 \includegraphics[scale=0.22]{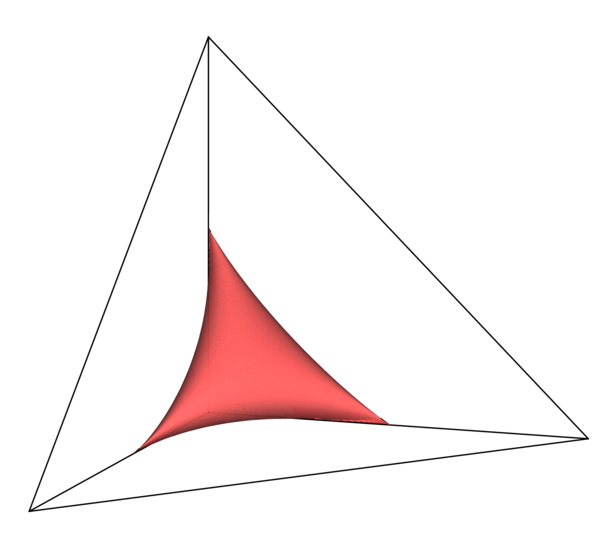} \\
 (a) $\mathbf{p}_1(\mathbf{s})$ & (b) $\mathbf{p}_2(\mathbf{s})$
 \end{tabular}
\caption{Examples of B\'ezier triangles $\mathbf{p}_1(\mathbf{s})$ defined in Section \ref{subsec:bez} and  $\mathbf{p}_2(\mathbf{s})$ defined in Section \ref{sec:bez2}. Exact and approximate implicitizations of the latter surface are in Figure \ref{fig:first}}
\label{fig:beziertri}
\end{center}
\end{figure}

Figure \ref{fig:beziertri}(a) shows an example of a degenerate quadratic B\'ezier triangle $\mathbf{p}_1(\mathbf{s}),$ with Cartesian control points
\begin{eqnarray*}
&\mathbf{c}_{200} = (1,0,0),& \\
&\mathbf{c}_{110} = (0,0,0), \ \ \mathbf{c}_{101} = (0,0,0),& \\
&\mathbf{c}_{020} = (0,1,0), \ \ \mathbf{c}_{011} = (0,0,0), \ \ \mathbf{c}_{002} = (0,0,1).&
\end{eqnarray*}
In Section \ref{sec:examples} we will see three alternative quadratic implicit approximations of this surface. Notice that the lexicographical ordering here is given by reading the control points from left to right and top to bottom.

\section{Approximate Implicitization}

In this section we outline the approach to approximate implicitization presented in \cite{dokken_thesis,dokken_approx}, in the context of B\'ezier triangles. Both the original approach and the so-called weak approach will be described, closely following the procedure given in \cite{dokken_weak}. We will also look at a numerical approach to the algorithm in greater detail.

The exact implicitization of a degree $n$ B\'ezier triangle may require a degree as high as $n^2.$ It should be noted that if a degree high enough for an exact implicitization is chosen and the algorithm is executed using exact arithmetic, then these methods will be exact. Use of floating point arithmetic will result in small rounding errors.

We begin by stating the formal definition of approximate implicitization:

An algebraic surface defined by the points $\mathbf{u}\in\mathbb{R}^3$ such that $q(\mathbf{u}) = 0$ for some polynomial $q,$ approximates the parametric surface $\mathbf{p}(\mathbf{s})$ within a tolerance of $\epsilon$ if there exists a vector-valued function $\mathbf{g}(\mathbf{s})$ of unit length, and an error function $\delta(\mathbf{s})$ such that 
\begin{equation}
q(\mathbf{p}(\mathbf{s}) + \delta(\mathbf{s})\mathbf{g}(\mathbf{s})) = 0,\label{eq:affineerr}
\end{equation}
and
\begin{displaymath}
\max_{\mathbf{s}\in\Omega} |\delta(\mathbf{s})| < \epsilon.
\end{displaymath}

We do not attempt to find the functions $\mathbf{g(s)}$ and $\delta(\mathbf{s})$ directly. Instead, we aim to find a polynomial $q$ of chosen degree $m$ that minimizes the algebraic distance $|q(\mathbf{p}(\mathbf{s}))|$ between the parametric and implicit surfaces. Certainly, if $q(\mathbf{p(s)})=0,$ then we have an exact implicitization. In Section \ref{sec:accuracy} we will see that this approach is also justified for approximate implicitization.

The method we use to find the polynomial $q,$ both in the original and weak approaches is to find the coefficients $b_\mathbf{i}$ of $q$ when expressed in the Bernstein basis of chosen degree $m:$ 
\begin{equation}\label{eq:factor1}
q(\mathbf{u}) = \sum_{|\mathbf{i}|=m} b_\mathbf{i} B_{\mathbf{i}}^{m}(\mathbf{u}).
\end{equation}

The difference between original and weak approximate implicitization is the choice of how to minimize the algebraic distance. The original approach attempts to minimize the pointwise error 
\begin{displaymath}
\max_{s\in\Omega} |q(\mathbf{p(s)})|,
\end{displaymath} 
whereas the weak approach minimizes by integration:
\begin{displaymath}
\int_{\Omega}(q(\mathbf{p(s)}))^2\ \mathrm{d}\mathbf{s}.
\end{displaymath}

\subsection{The Original Approach}
We follow the same steps as in approximate implicitization of tensor-product B\'ezier surfaces and B\'ezier curves, only now using the triangular Bernstein basis functions. As we will see, the details differ somewhat in the triangular case. 

Since we have chosen $q$ to be of degree $m,$ and $\mathbf{p(s)}$ is defined to be degree $n,$ the expression $q(\mathbf{p(s)})$ will be a polynomial of degree $mn$ in $\mathbf{s}.$  Such a polynomial can be factorized in the Bernstein basis $(B_{\mathbf{j}}^{mn})_{|\mathbf{j}|=mn}$ with coefficients $d_\mathbf{i,j}.$ To obtain these coefficients, we form the following composition of the coordinate functions of $\mathbf{p}(\mathbf{s})$ with each Bernstein basis function $B_\mathbf{i}^m:$
\begin{equation}\label{eq:dcolumns}
B_\mathbf{i}^{m}(\mathbf{p(s)}) = \sum_{|\mathbf{j}|=mn} d_\mathbf{i,j} B_\mathbf{j}^{mn}(\mathbf{s}).
\end{equation}
Note that $ d_\mathbf{i,j}$ can be calculated explicitly by using (\ref{eq:product}), the product rule for Bernstein bases. An example of how this is done is presented in Section \ref{subsec:firstex}.

Now, using (\ref{eq:factor1}) and (\ref{eq:dcolumns}) we get
\begin{eqnarray}
q(\mathbf{p(s}))& = & \sum_{|\mathbf{i}|=m} b_\mathbf{i} B_{\mathbf{i}}^{m}(\mathbf{p(s}))\nonumber \\
& = & \sum_{|\mathbf{i}|=m} b_\mathbf{i} \left( \sum_{|\mathbf{j}|=mn} d_\mathbf{i,j} B_\mathbf{j}^{mn}(\mathbf{s}) \right) \nonumber \\ 
& = & \sum_{|\mathbf{j}|=mn} B_\mathbf{j}^{mn}(\mathbf{s}) \left( \sum_{|\mathbf{i}|=m} d_\mathbf{i,j} b_\mathbf{i} \right). \label{eq:factor4}
\end{eqnarray}
Since the matrix $\mathbf{D}$ defined by the coefficients $(d_\mathbf{i,j})_{|\mathbf{i}|=m,|\mathbf{j}|=mn}$ is fundamental to the theory of approximate implicitization, we summarize its construction in the following proposition:

\begin{prop}
The ${ \binom{m + 3}{3} } \times { \binom{mn+2}{2} }$ matrix $\mathbf{D}$ for approximate implicitization of triangular B\'ezier surfaces can be constructed by repeated multiplication of the coordinate functions of $\mathbf{p}(\mathbf{s}),$ according to the equation (\ref{eq:dcolumns}).
\end{prop}


Writing the unknown coefficients $(b_\mathbf{i})_{|\mathbf{i}|=m}$ and the basis functions $(B_{\mathbf{j}}^{mn}(\mathbf{s}))_{|\mathbf{j}|=mn}$ in vectors $\mathbf{b}$ and $\mathbf{B}^{mn}(\mathbf{s})$ respectively, we restate (\ref{eq:factor4}) as
\begin{equation}
q(\mathbf{p(s})) = \mathbf{B}^{mn}(\mathbf{s})^T\mathbf{Db}.\label{eq:factor3}
\end{equation}
We may impose, without loss of generality, the normalization condition $\Vert\mathbf{b}\Vert=1.$ Since the Bernstein basis forms a partition of unity, using the factorization (\ref{eq:factor3}) we get
\begin{eqnarray}
\max_{\mathbf{s}\in\Omega}|q(\mathbf{p(s)})| & = & \max_{\mathbf{s}\in\Omega}|\mathbf{B}^{mn}(\mathbf{s})^T\mathbf{Db}| \nonumber \\
& \leq & \max_{\mathbf{s}\in\Omega}\Vert \mathbf{B}^{mn}(\mathbf{s})\Vert\Vert\mathbf{Db}\Vert\leq\Vert\mathbf{Db}\Vert. \nonumber
\end{eqnarray}
The approximation may well be good outside the region of interest $\Omega,$ but the result used here, that $\Vert \mathbf{B}^{mn}(\mathbf{s})\Vert\leq 1,$ is specific only to the domain $\Omega.$
A standard result from linear algebra tells us that $\min_{\Vert \mathbf{b} \Vert = 1}\Vert \mathbf{Db} \Vert = \sigma_{\min},$ where $\sigma_{\min}$ is the smallest singular value of $\mathbf{D}.$ So, in particular we have
\begin{equation}
\min_{\Vert\mathbf{b}\Vert=1}\max_{\mathbf{s}\in\Omega}|q(\mathbf{p(s)})|\leq \sigma_{\min}.
\end{equation}
We can thus minimize the left hand side of the inequality by performing a singular value decomposition (SVD) on the matrix $\mathbf{D}.$ The vector $\mathbf{b}_{\min}$ corresponding to the smallest singular value $\sigma_{\min}$ of $\mathbf{D}$ would then give the best candidate for the approximation.

\subsection{The Weak Approach}\label{sec:weak}

Recall that the weak approach attempts to minimize the algebraic distance by minimizing the integral $\int_{\Omega}(q(\mathbf{p(s)}))^2\ \mathrm{d}\mathbf{s}.$ Here we approach this problem using the exact integration formula (\ref{eq:basisint}). However, the weak approach also introduces the possibility to perform a numerical integration. In Section \ref{sec:numerical} we will discuss this further.

Using the factorization (\ref{eq:factor3}) we can perform the integral as follows:
\begin{eqnarray}
\int_{\Omega} \left(q(\mathbf{p(s)})\right)^2 \ \mathrm{d}\mathbf{s} & = & \int_{\Omega}  \left(\mathbf{B}^{mn}(\mathbf{s})^T\mathbf{Db} \right)^2 \ \mathrm{d}\mathbf{s} \nonumber \\
& = & \mathbf{b}^T \mathbf{D}^T \left(\int_{\Omega} \mathbf{B}^{mn}(\mathbf{s})^T\mathbf{B}^{mn}(\mathbf{s}) \ \mathrm{d}\mathbf{s}\right)  \mathbf{Db} \nonumber \\
& = & \mathbf{b}^T \mathbf{D}^T \mathbf{A} \mathbf{Db},
\end{eqnarray} 
where $\mathbf{A}$ is the symmetric matrix defined by $(a_{\mathbf{i,j}})_{|\mathbf{i}|=mn,|\mathbf{j}|=mn}$
\begin{eqnarray*}
a_{\mathbf{i},\mathbf{j}} & = & \int_\Omega B_{\mathbf{i}}^{mn}(\mathbf{s}) B_{\mathbf{j}}^{mn}(\mathbf{s}) \ \mathrm{d}\mathbf{s} \\
& = & \frac{{ \binom{mn}{\mathbf{i}}}{ \binom{mn}{\mathbf{j}}}}{ \binom{2mn}{\mathbf{i+j}}}  \int_\Omega B_{\mathbf{i+j}}^{2mn}(\mathbf{s}) \ \mathrm{d}\mathbf{s} \\
& = & \frac{{ \binom{mn}{\mathbf{i}} }{ \binom{mn}{\mathbf{j}} }}{{ \binom{2mn}{\mathbf{i+j}} }} \frac{1}{(2mn+1)(2mn+2)}.
\end{eqnarray*}
We may define the matrix $\mathbf{M}$ by 
\begin{equation}
\mathbf{M} = \mathbf{D}^T \mathbf{A} \mathbf{D}.\label{eq:DAD}
\end{equation}
Then, similarly to the original approach, an SVD of $\mathbf{M}$ will give us a candidate for a weak approximate implicitization of $\mathbf{p(s)}.$ We again choose the vector corresponding to the smallest singular value for the best candidate. The construction of $\mathbf{M}$ is summarized as follows: 
\begin{prop}
The  ${ \binom{m+3}{3} } \times { \binom{m+3}{3} }$ matrix $\mathbf{M}$ formed in weak approximate implicitization of triangular B\'ezier surfaces can be built by the matrix multiplication $\mathbf{D}^T\mathbf{A}\mathbf{D},$ where the matrix $\mathbf{A}$ depends only on $m$ and $n.$
\end{prop}
Since $\mathbf{A}$ is only dependent on the degrees $m$ and $n,$ it could in fact be pre-calculated, meaning the construction of $\mathbf{M}$ is reduced to making two matrix multiplications. 

This method may be particularly useful when combining the original and weak approximations in order to remove unwanted branches, as the $\mathbf{D}$ matrix must already be calculated. By combining the best approximations from the original and weak forms, we will obtain another approximation with a high convergence rate. Since both the approximations will be `good' in the area of interest, but may have different branches, the combination may remove these unwanted branches. 

For a detailed discussion of the relationship between the weak and original forms of approximate implicitization, we refer the reader to \cite{dokken_weak}. Here we simply state the main results:
\begin{displaymath}
|q(\mathbf{p(s}))| \leq \frac{1}{\sqrt{\lambda_{\min}}} \Vert \bm{\Sigma} \mathbf{UDb}\Vert,
\end{displaymath}
and 
\begin{displaymath}
\sqrt{\int_\Omega \left(q(\mathbf{p(s}))\right)^2 \ \mathrm{d}\mathbf{s}} \leq \sqrt{\lambda_{\max}} \Vert \mathbf{Db} \Vert,
\end{displaymath}
where $\bm{\Sigma}$ is a diagonal matrix containing the square roots of the eigenvalues $\lambda_{\min}, \ldots, \lambda_{\max}$ of $\mathbf{A},$ and $\mathbf{A}=\mathbf{U}^T(\bm{\Sigma}^2)\mathbf{U}.$

\subsection{Numerical Approximation}\label{sec:numerical}

As the exact integration in weak approximate implicitization can be replaced by a numerical integration, the need for an explicit rational parametric form is removed. Numerical integration only requires that the surface can be evaluated. This allows, for example, procedural surfaces to be approximated. Integration using numerical methods allows for quick building of the $\mathbf{M}$ matrix. In addition, we show that the algorithm exhibits symmetries that further enhance its efficiency. The results of this section can be easily generalized to apply to weak approximate implicitization of rational parametric manifolds in $\mathbb{R}^l.$

In the previous section we constructed $\mathbf{M}$ via matrix multiplications. Perhaps a more natural method to construct $\mathbf{M}$ is to perform the integration using the equation (\ref{eq:factor1}). Using this method we obtain an element-wise formula for $\mathbf{M},$ which we can evaluate by making use of (\ref{eq:product}):
\begin{eqnarray}
m_{\mathbf{i,j}} & = & \int_{\Omega} B_{\mathbf{i}}^m(\mathbf{p(s)})B_{\mathbf{j}}^m(\mathbf{p(s)}) \ \mathrm{d}\mathbf{s} \label{eq:mij}\\
& = & \frac{{ \binom{m}{\mathbf{i}}}{ \binom{m}{\mathbf{j}}}}{ \binom{2m}{\mathbf{i+j}}} \int_{\Omega} B_{\mathbf{i+j}}^{2m}(\mathbf{p(s)}) \ \mathrm{d}\mathbf{s}. \label{eq:appint1}
\end{eqnarray}
This method eliminates the need to compute $\mathbf{D},$ but we are now required to evaluate the polynomials $B_{\mathbf{i+j}}^{2m}(\mathbf{p(s)})$ in order to use (\ref{eq:bernint}) for the integration. This is in comparison to evaluating the expressions $B_{\mathbf{i}}^m(\mathbf{p(s)})$ required to build $\mathbf{D}.$ Due to the lower polynomial degrees involved in the latter, the construction of $\mathbf{M}$ by first computing $\mathbf{D}$ and then applying (\ref{eq:DAD}), is preferable for the exact integration. However, equation (\ref{eq:appint1}) provides a direct method that would be preferable if using numerical integration, since it avoids the polynomial multiplication. 

Inspecting (\ref{eq:mij}) we clearly see that $\mathbf{M}$ is symmetric. However, there exist other symmetries which allow us to avoid repeated calculation of the integrals for each element $m_{\mathbf{i,j}}.$ Equation (\ref{eq:appint1}) shows that there are in fact only $\binom{2m+3}{3}$ unique integrals required. We can thus pre-calculate these integrals using some chosen numerical integration method:  
\begin{equation}
\left(\int_{\Omega} B_{\mathbf{k}}^{2m}(\mathbf{p(s)})\right)_{|\mathbf{k}|=2m}.\label{eq:integrals}
\end{equation}
Exploiting these symmetries results in the required number of integrals being proportional to $m^3$ rather than $m^6.$ This result is summarized in the following proposition:
\begin{prop}
The $\binom{m+3}{3} \times \binom{m+3}{3}$ matrix $\mathbf{M}$ formed in weak approximate implicitization of triangular B\'ezier surfaces can be built by pre-computing the $\binom{2m+3}{3}$ integrals in (\ref{eq:integrals}), and multiplying the relevant integrals with the coefficients $\frac{ \binom{m}{\mathbf{i}} \binom{m}{\mathbf{j}} }{ \binom{2m}{\mathbf{i+j}} }.$
\end{prop}
Since the degree of the integrand is $2mn,$ it is vital that the numerical integration techniques used, exhibit numerical stability up to high polynomial degrees. For example, approximating a cubic B\'ezier triangle by a cubic implicit surface, requires the numerical integration of a bivariate polynomial of degree 18.



\subsection{Approximating Rational B\'ezier Triangles}

Rational B\'ezier triangles give extra flexibility in CAGD and are in fact required to be able to represent general quadric surfaces exactly. Although this can be done with rational tensor-product patches, some degeneracy is necessary, and hence singularities are introduced. We will show in this section that the algorithm for approximate implicitization of rational B\'ezier triangles is only a short extension of the non-rational version. We first introduce the concept of rational B\'ezier triangles, as described in \cite{farin_cagd}.

A rational B\'ezier triangle of degree $n$ is defined similarly to the non-rational case as follows:
\begin{displaymath}                      
\mathbf{r}(\mathbf{s}) = \sum_{|\mathbf{i}|=n} {\mathbf{c}_\mathbf{i} R_\mathbf{i}^n(\mathbf{s})},
\end{displaymath}
where 
\begin{displaymath}                      
R_{\mathbf{i}}^n(\mathbf{s}) = \frac{{w_\mathbf{i} B_\mathbf{i}^n(\mathbf{s})}}{\sum_{|\mathbf{i}|=n}{w_\mathbf{i}B_\mathbf{i}^n(\mathbf{s})}} = \frac{g_\mathbf{i}(\mathbf{s})}{h(\mathbf{s})}.
\end{displaymath}
The $w_{\mathbf{i}}$ denote weights assigned to each control point $\mathbf{c}_{\mathbf{i}}.$ Note that the basis $(R_\mathbf{i}(\mathbf{s}))_{|\mathbf{i}|=n}$ defines a partition of unity, so the original approach to approximate implicitization can be used in a similar way for rational B\'ezier triangles. In fact, on forming the expression $q(\mathbf{r}(\mathbf{s})),$ we can factor out the denominator, which allows us to consider only the numerator for an exact implicitization \cite{thomassen_selfint}. Since the numerator is simply a regular B\'ezier triangle (albeit with the weights absorbed into the control points), this implies that we can find implicitly defined quadrics from non-rational B\'ezier triangles. We show this as follows:
\begin{eqnarray} 
q(\mathbf r(\mathbf s)) & = & \sum_{|\mathbf j|=m} b_{\mathbf j} B_{\mathbf{j}}^m (\mathbf{r(s)}) \nonumber \\
& = & \frac{1}{(h(\mathbf{s}))^m} \sum_{|\mathbf j|=m} b_{\mathbf j} B_{\mathbf{j}}^m \left(\sum_{|\mathbf{i}|=n} {\mathbf{c}_\mathbf{i} g_\mathbf{i}(\mathbf{s})}\right). \label{eq:denominator}
\end{eqnarray}
We obtain an exact implicitization if and only if the sum over $|\mathbf{i}|=n$ in (\ref{eq:denominator}) is zero; but these are exactly the same conditions for exact implicitization on non-rational B\'ezier triangles, allowing us to disregard $h(\mathbf{s})$ in the algorithm. That is, we may perform the implicitization on
\begin{displaymath}
\sum_{|\mathbf i| = n} \mathbf{c}_{\mathbf i} g_{\mathbf{i}}(\mathbf{s}) =  \sum_{|\mathbf i| = n} \bm{\gamma}_{\mathbf i} B_{\mathbf{i}}^n(\mathbf{s}),
\end{displaymath}
where $\bm{\gamma}_\mathbf{i} = w_\mathbf{i}\mathbf{c}_\mathbf{i}.$
We may also disregard $h$ for approximate implicitizations, however this will come at some expense to the quality of approximation if the function $h$ has large variations.


\subsection{Accuracy in Affine Space and Convergence Rates}\label{sec:accuracy}

The intention of this section is to show why approximate implicitization works, and to state a result about the quality of the approximation. For a more in-depth coverage of these topics see \cite{dokken_thesis}.

Recall the definition of approximate implicitization from the beginning of this section. This definition ensures that the implicit and parametric curves lie close together in affine space. However, by minimizing the algebraic distance, as we did in the algorithm, we cannot necessarily guarantee that the affine error will be small. The affine and algebraic errors are related by the following Taylor expansion of (\ref{eq:affineerr}):
\begin{displaymath}
q(\mathbf{p}(\mathbf{s})) + \delta(\mathbf{s})\mathbf{g}(\mathbf{s}) \cdot \nabla q(\mathbf{p}(\mathbf{s}))  + \cdots = 0.
\end{displaymath}
Suppose we have a polynomial $q$ such that $q(\mathbf{p}(\mathbf{s})) \approx 0.$ Then the above equation shows that either  $\nabla q(\mathbf{p}(\mathbf{s}))$ or $\delta(\mathbf{s})$ must be small. Certainly, away from singularities, where the gradient $\nabla q(\mathbf{p}(\mathbf{s}))$ does not vanish, $\delta(\mathbf{s})$ will be small, meaning the approximation in affine space is good. This justifies the approach to approximate implicitization outlined above, away from singularities. In the region of singularities, the neighbourhood of the singular point or curve will attract the approximation to the correct shape; however, the singularities themselves may be smoothed out. A clear example of this is the approximation in Figure \ref{fig:first}. Here, the approximation is attracted to the non-singular part of the surface and the singular curves are `smoothed out'. We will consider this example further in Section \ref{sec:bez2}.

We can improve the approximation in affine space by performing the approximation over a smaller region of the parametric surface. The convergence rates of approximate implicitization, as the size of the region to be approximated is reduced, have been investigated in \cite{dokken_overview}. Here we state the result most relevant to this paper; the convergence rate of surfaces in $\mathbb{R}^3.$ Given a closed box of diameter $h$ in $\Omega$ around a point $\mathbf{s}_0,$ we have the convergence rate
\begin{equation}
O\left(h^{\lfloor \frac{1}{6} \sqrt{(9+12m^3+72m^2+132m)} \rfloor-\frac{1}{2}}\right).
\end{equation}
Here, $\lfloor x \rfloor$ denotes the integer part of $x.$

\section{Examples of Implicitization of B\'ezier Triangles}\label{sec:examples}

In this section we present examples of approximate implicitization of triangular B\'ezier surfaces. We begin with a simple example that can be calculated by hand, before moving on to more computationally intensive examples. Our first example will find an implicit surface that approximates a single quadratic B\'ezier triangle.

\subsection{A First Example}\label{subsec:firstex}

Recall the definition of $\mathbf{p}_1(\mathbf{s}),$ the degenerate quadratic B\'ezier surface mentioned in Section \ref{subsec:bez} and pictured in Figure \ref{fig:beziertri}(a). The control points also form a tetrahedron over which we can define the barycentric coordinate system. Using these barycentric coordinates, the patch is described by,
\begin{equation}
\mathbf{p}_1(\mathbf{s}) = (B_{200}^2(\mathbf{s}),B_{020}^2(\mathbf{s}),B_{002}^2(\mathbf{s}),B_{011}^2(\mathbf{s})+B_{101}^2(\mathbf{s})+B_{110}^2(\mathbf{s})). \label{eq:bezier1}
\end{equation}
For this example,  we choose to approximate $\mathbf{p}_1(\mathbf{s})$ by a quadratic implicit surface, in order to keep the matrix $\mathbf{D}$ to a manageable size. However, an exact implicitization in fact requires an implicit surface of degree four.

A trivariate polynomial of degree two, represented in the tetrahedral Bernstein basis can be written as follows:
\begin{displaymath}
q(\mathbf{u}) = \sum_{|\mathbf{i}|=2} b_\mathbf{i} B_{\mathbf{i}}^{2}(\mathbf{u}),
\end{displaymath}
for barycentric coordinates $\mathbf{u}.$

Now, to construct the $15 \times 10$ matrix $\mathbf{D},$ we simply expand the expression (\ref{eq:dcolumns}), for each of the basis functions $(B_\mathbf{i}^2(\mathbf{u}))_{|\mathbf{i}|=2}$ and write the resulting coefficients in the columns of $\mathbf{D}$. We use the lexicographical ordering system to relate the entries of the matrix, to the multi-indices.

The first column in the matrix $\mathbf{D}$ contains the coefficients of $B_{2000}^2(\mathbf{p}_1(\mathbf{s)})$ which by the definition (\ref{eq:bezier1}) and the product rule (\ref{eq:product}) is equal to $(B_{200}^2(\mathbf{s}))^2=B_{400}^4(\mathbf{s}).$ The first column of $\mathbf{D}$ is thus the vector of coefficients that are all zero except for the first:
\[
(1,0,0,0,0,0,0,0,0,0,0,0,0,0,0)^T.
\] 
Similarly, the second column is calculated by expanding $B_{1100}^2(\mathbf{p}_1(\mathbf{s)})$ giving 
\[
(0,0,0,\frac{1}{3},0,0,0,0,0,0,0,0,0,0,0)^T.
\]
Continuing in this way we get the matrix,
\begin{displaymath}
\mathbf{D} = \left( \begin{matrix}
1 & 0 & 0 & 0 & 0 & 0 & 0 & 0 & 0 & 0 \\
0 & 0 & 0 & 1 & 0 & 0 & 0 & 0 & 0 & 0 \\
0 & 0 & 0 & 1 & 0 & 0 & 0 & 0 & 0 & 0 \\
0 & \frac{1}{3} & 0 & 0 & 0 & 0 & 0 & 0 & 0 & \frac{2}{3} \\
0 & 0 & \frac{1}{3} & 0 & 0 & 0 & 0 & 0 & 0 & \frac{2}{3} \\
0 & 0 & 0 & \frac{1}{3} & 0 & 0 & 0 & 0 & 0 & \frac{2}{3} \\
0 & 0 & 0 & 0 & 0 & 0 & 1 & 0 & 0 & 0 \\
0 & 0 & 0 & 0 & 0 & 0 & 0 & 0 & 1 & 0 \\
0 & 0 & 0 & 0 & 0 & 0 & \frac{1}{3} & 0 & 0 & \frac{2}{3} \\
0 & 0 & 0 & 0 & 0 & 0 & 0 & 0 & \frac{1}{3} & \frac{2}{3} \\
0 & 0 & 0 & 0 & 1 & 0 & 0 & 0 & 0 & 0 \\
0 & 0 & 0 & 0 & 0 & 0 & 1 & 0 & 0 & 0 \\
0 & 0 & 0 & 0 & 0 & \frac{1}{3} & 0 & 0 & 0 & \frac{2}{3} \\
0 & 0 & 0 & 0 & 0 & 0 & 0 & 0 & 1 & 0 \\
0 & 0 & 0 & 0 & 0 & 0 & 0 & 1 & 0 & 0 
\end{matrix} \right). \nonumber
\end{displaymath}
The correct and accurate construction of this matrix can be confirmed by checking that the rows sum to 1 (see Theorem 4.3 in \cite{dokken_approx}). As we have proceeded using exact methods, we expect no errors here.

\begin{figure}  
\begin{center}
  \includegraphics[width=1.7in]{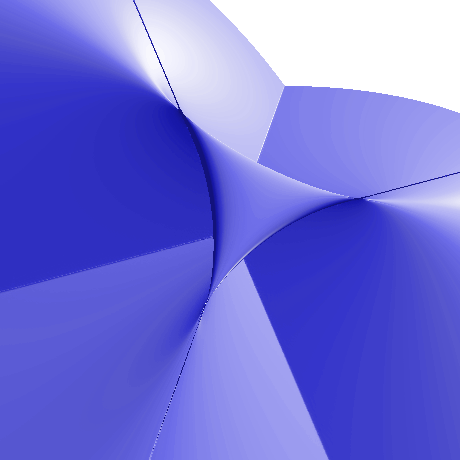}\hspace*{1cm}
  \includegraphics[width=1.7in]{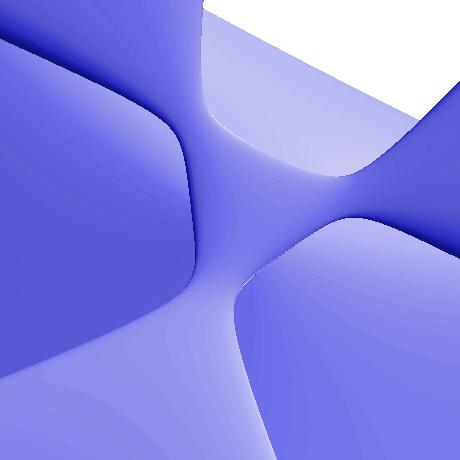}
  \caption[Exact and approximate implicitizations of a B\'ezier triangle]{Exact (left) and approximate (right) implicitizations of a quadratic B\'ezier triangle with 
  singularities $\mathbf{p}_2(\mathbf{s})$ (see Figure \ref{fig:beziertri}(b)).}\label{fig:first}
\end{center}
\end{figure}

We now perform an SVD on this matrix, and choose the vector $\mathbf{b}$ corresponding to the smallest singular value $\sigma_{\min}.$ The singular values of $\mathbf{D}$ are
\begin{displaymath}
\begin{split}
(1.70471, 1.45296, 1.45296, 1.38925, 1.00000,\\ 1.00000, 1.00000, 0.33333, 0.33333, 0.&22984),
\end{split}
\end{displaymath}
and the normalized vector corresponding to $\sigma_{\min} = 0.22984$ is 
\begin{displaymath}
\begin{split}
\mathbf{b}_{\textrm{orig}} = (0.00000, -0.57062, -0.57062, -0.01616, 0.00000, \\-0.57062, -0.01616, 0.00000, -0.01616, 0.1496&6).
\end{split}
\end{displaymath}
This vector defines a candidate for an approximate implicitization of $\mathbf{p}_1(\mathbf{s}).$ 

The approach of weak approximate implicitization is equally well suited to this example. As stated previously, property (\ref{eq:bernint}) allows us to integrate Bernstein polynomials by summing the coefficients of the Bernstein basis and dividing by a constant factor. For simplicity, we proceed here using the element-wise definition of $\mathbf{M},$ (\ref{eq:appint1}):
\begin{eqnarray*}
m_\mathbf{i,j} & = & \frac{ \binom{2}{\mathbf{i}} \binom{2}{\mathbf{j}} }{ \binom{4}{\mathbf{i+j}} } \int_{\Omega} B_{\mathbf{i+j}}^{4}(\mathbf{p}_1(\mathbf{s})) \ \mathrm{d}\mathbf{s}. \\
\end{eqnarray*}
For example, the first entry for $\mathbf{i}=\mathbf{j}=(2,0,0,0),$ is calculated by making the integration 
\begin{displaymath}
\int_{\Omega} B_{4000}^{4}(\mathbf{\mathbf{p}_1(\mathbf{s})}) = \int_{\Omega} (B_{200})^4 = \int_{\Omega} B_{800}^8.
\end{displaymath}
This is a degree eight Bernstein polynomial with first coefficient equal to one and all other coefficients equal to zero. The first value of the matrix is thus $m_{1,1} = 1/90.$ The other values of the matrix can be computed similarly to get
\begin{displaymath}
\mathbf{M} = \left( \begin{smallmatrix}
\frac{1}{90} & \frac{1}{1260} & \frac{1}{1260} & \frac{1}{84} & \frac{1}{6300} & \frac{1}{18900} & \frac{23}{18900} & \frac{1}{6300} & \frac{23}{18900} & \frac{4}{675}\cr \frac{1}{1260} & \frac{1}{1575} & \frac{1}{9450} & \frac{23}{9450} & \frac{1}{1260} & \frac{1}{9450} & \frac{23}{9450} & \frac{1}{18900} & \frac{1}{2100} & \frac{31}{9450}\cr \frac{1}{1260} & \frac{1}{9450} & \frac{1}{1575} & \frac{23}{9450} & \frac{1}{18900} & \frac{1}{9450} & \frac{1}{2100} & \frac{1}{1260} & \frac{23}{9450} & \frac{31}{9450}\cr \frac{1}{84} & \frac{23}{9450} & \frac{23}{9450} & \frac{16}{675} & \frac{23}{18900} & \frac{1}{2100} & \frac{31}{4725} & \frac{23}{18900} & \frac{31}{4725} & \frac{67}{3150}\cr \frac{1}{6300} & \frac{1}{1260} & \frac{1}{18900} & \frac{23}{18900} & \frac{1}{90} & \frac{1}{1260} & \frac{1}{84} & \frac{1}{6300} & \frac{23}{18900} & \frac{4}{675}\cr \frac{1}{18900} & \frac{1}{9450} & \frac{1}{9450} & \frac{1}{2100} & \frac{1}{1260} & \frac{1}{1575} & \frac{23}{9450} & \frac{1}{1260} & \frac{23}{9450} & \frac{31}{9450}\cr \frac{23}{18900} & \frac{23}{9450} & \frac{1}{2100} & \frac{31}{4725} & \frac{1}{84} & \frac{23}{9450} & \frac{16}{675} & \frac{23}{18900} & \frac{31}{4725} & \frac{67}{3150}\cr \frac{1}{6300} & \frac{1}{18900} & \frac{1}{1260} & \frac{23}{18900} & \frac{1}{6300} & \frac{1}{1260} & \frac{23}{18900} & \frac{1}{90} & \frac{1}{84} & \frac{4}{675}\cr \frac{23}{18900} & \frac{1}{2100} & \frac{23}{9450} & \frac{31}{4725} & \frac{23}{18900} & \frac{23}{9450} & \frac{31}{4725} & \frac{1}{84} & \frac{16}{675} & \frac{67}{3150}\cr \frac{4}{675} & \frac{31}{9450} & \frac{31}{9450} & \frac{67}{3150} & \frac{4}{675} & \frac{31}{9450} & \frac{67}{3150} & \frac{4}{675} & \frac{67}{3150} & \frac{22}{525}
\end{smallmatrix} \right). \nonumber
\end{displaymath}
The accuracy of the construction of this matrix can be confirmed by checking that the elements sum to $\frac{1}{2}$ (see Theorem 2 in \cite{dokken_weak}).

\begin{figure*}
\begin{center}
\hspace*{0.8cm}
\begin{tabular}{cccc}
\begin{minipage}[b]{0.25\linewidth}
	\centering
	\includegraphics[scale=0.2]{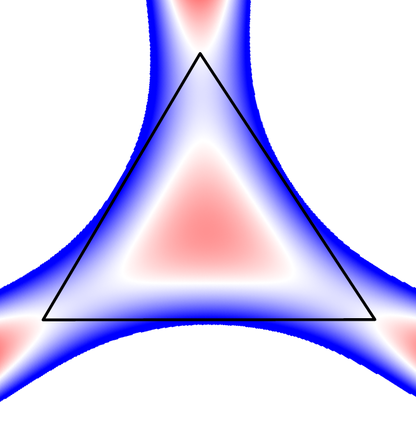}{\begin{center}Original\end{center}}\label{fig:errors1}
\end{minipage}			
\begin{minipage}[b]{0.25\linewidth}
	\centering
	 \includegraphics[scale=0.2]{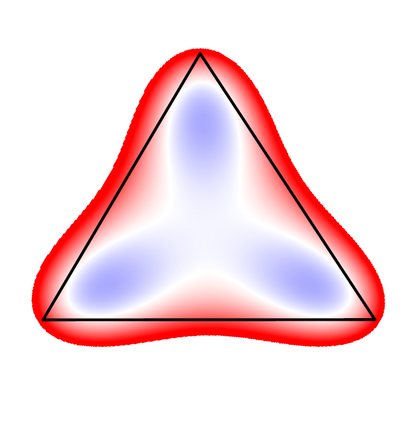}{\begin{center}Weak\end{center}}\label{fig:errors2}
\end{minipage}		
\begin{minipage}[b]{0.25\linewidth}
	\centering
	 \includegraphics[scale=0.2]{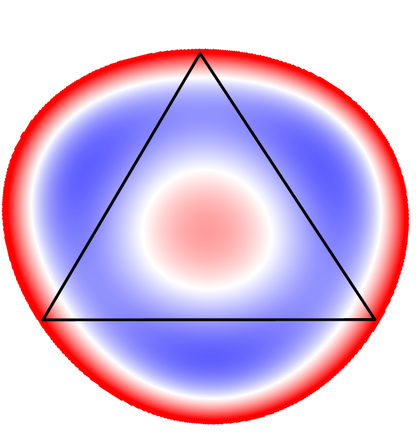}{\begin{center}Combined\end{center}}\label{fig:errors3}
\end{minipage}
\hspace*{-1.0cm}												
\begin{minipage}[b]{0.25\linewidth}
	\centering
	\includegraphics[scale=0.13]{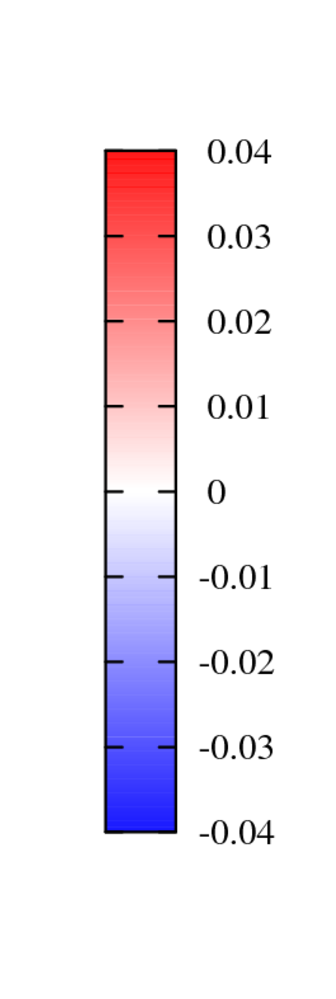}{ }
\end{minipage}				
\end{tabular}
\caption[Colour maps of the algebraic error fo the approximations]{Colour maps showing the algebraic approximation errors $q_{\mathbf{b}}(\mathbf{p}_1(\mathbf{s}))$ 
over the domain $\Omega$ bounded by the triangle. Note that the white parts correspond to intersection curves between the implicit and parametric surfaces.}\label{fig:errors}
\end{center}
\end{figure*}

Again, performing an SVD on this matrix and choosing the vector corresponding to the smallest singular value will define an implicit equation that is a candidate for approximation:
\begin{displaymath}
\begin{split}
\mathbf{b}_{\textrm{weak}} = (0.03985, 0.56837, 0.56837, -0.09313, 0.03985,\\ 0.56837, -0.09313, 0.03985, -0.09313, -0.00&859).
\end{split}
\end{displaymath}

Although this simple example has no extraneous branches, in order to illustrate the possibility of modelling the shape of the approximation, we include a combined approximation. This is obtained by summing the coefficients of the original and weak approximations and renormalizing:
\begin{displaymath}
\begin{split}
\mathbf{b}_{\textrm{comb}} = (0.11496, -0.00652, -0.00652, -0.31523, 0.11496,\\ -0.00652, -0.31523, 0.11496, -0.31523, 0.8137&1).
\end{split}
\end{displaymath}

Figure \ref{fig:errors} shows the algebraic distance between the parametric and approximate implicit surfaces. The three approximations exhibit different behaviour with regard to where the surfaces intersect and the positions of the maximum error. This illustrates the possibility of modelling the surfaces to obtain certain characteristics. Alternative approximations could also be formed by taking different combinations of the two surfaces, by combining approximations corresponding to other singular values in the SVD, or by adding constraints to the algebraic equation.

When constructing this example we ensured that the corners of the B\'ezier triangle were reused as vertices in the tetrahedral barycentric coordinate system, with the remaining fourth vertex positioned symmetrically with respect to these three points. This symmetry is reflected in the intersection curves between the triangular B\'ezier surfaces and the approximations in Figure \ref{fig:errors}. In this example, the original approximate implicitization intersects the corners of the triangular B\'ezier surface; however, the interpolation is special to this case. It is easy to construct examples with the same collocation of surface corners and tetrahedral vertices where the approximate implicit generated by the original approach does not intersect the corners of the triangular B\'ezier surface.

\subsection{Approximation of Several Patches with One Implicit Surface}\label{subsec:several}
\begin{figure*}[ht]
\begin{tabular}{cccc}
\begin{minipage}[b]{0.23\linewidth}
	\centering
	\includegraphics[width=3.0cm,height=2cm]{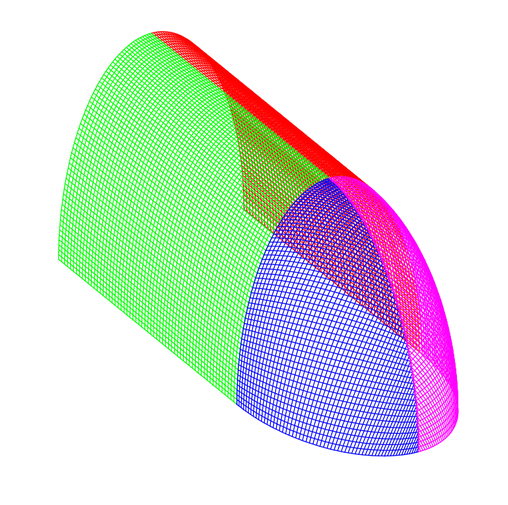}{\begin{center}(a) parametric form \phantom{parametric} \end{center}}\label{fig:parametric}
\end{minipage}		
\begin{minipage}[b]{0.23\linewidth}
	\centering
	\includegraphics[width=3.2cm]{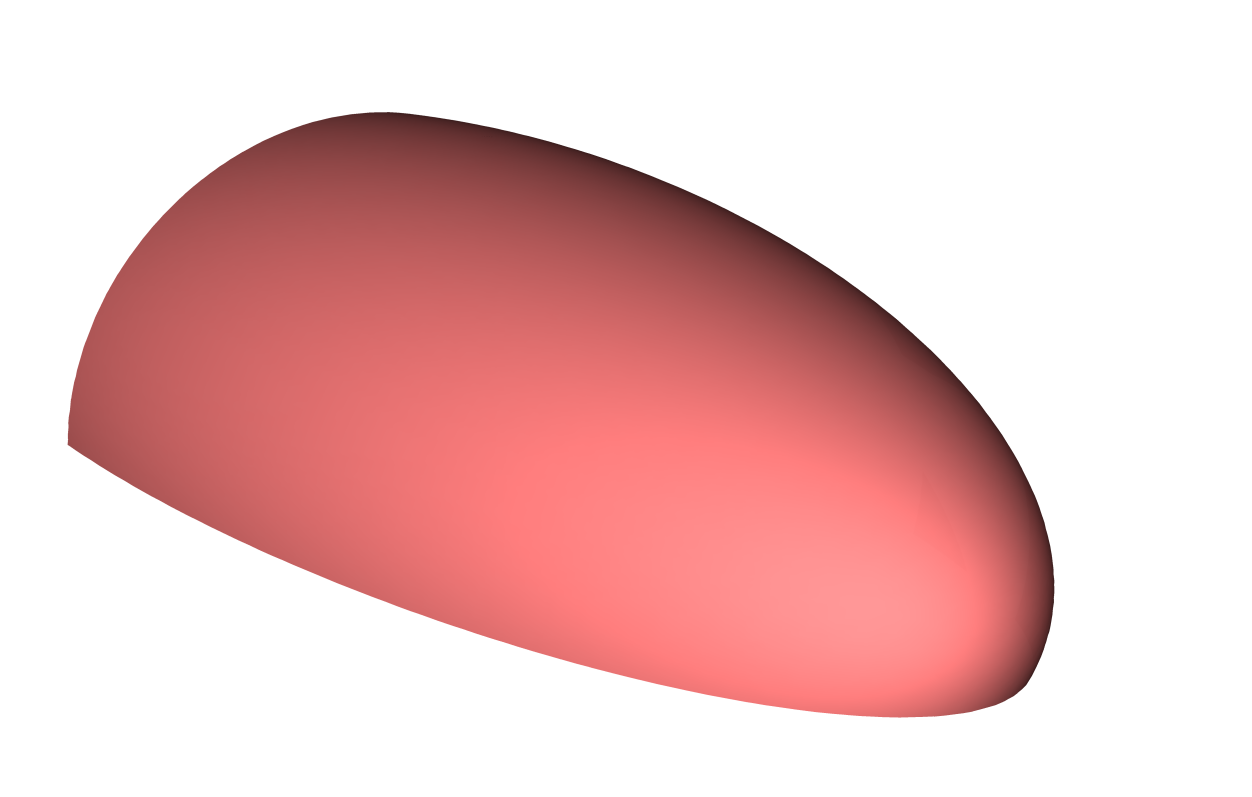}{\begin{center}(b) quadratic approximation\end{center}}\label{fig:quadratic}
\end{minipage}			
\begin{minipage}[b]{0.23\linewidth}
	\centering
    \includegraphics[width=3.2cm]{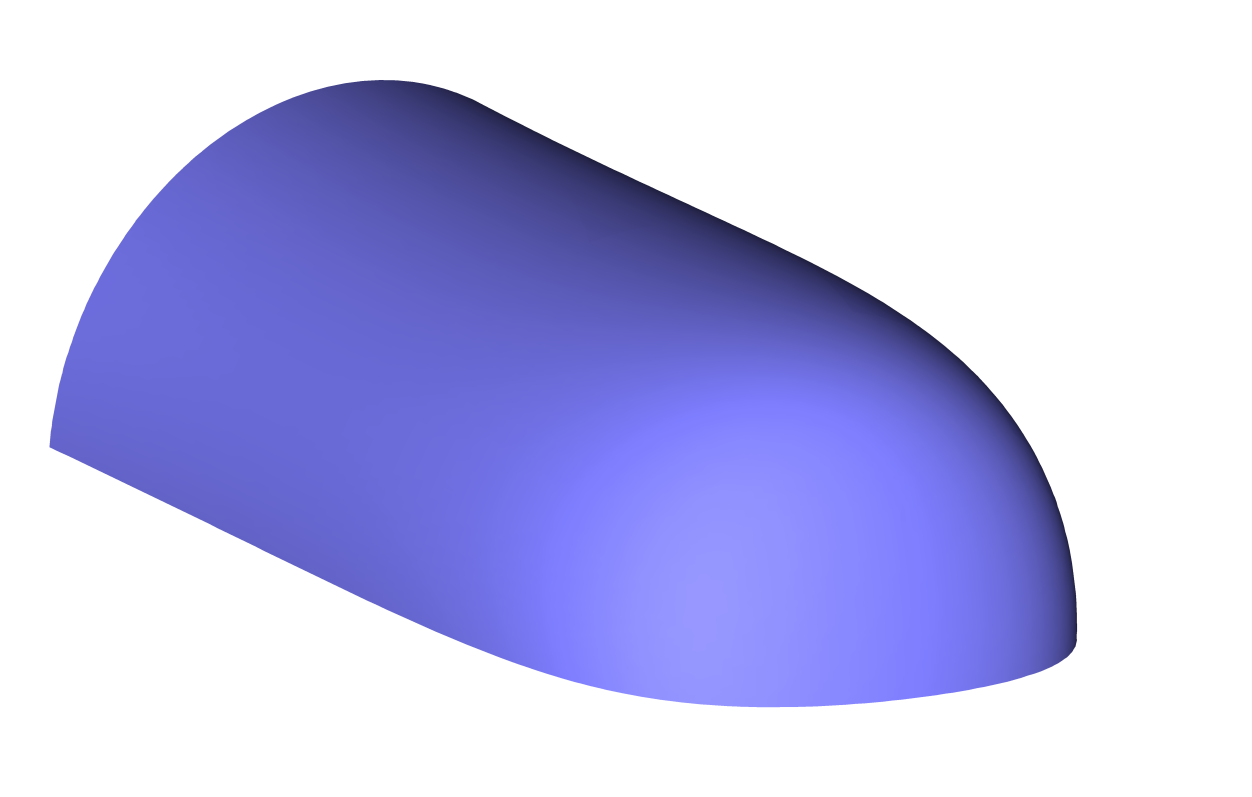}{\begin{center}(c) cubic approximation\end{center}}\label{fig:cubic}
\end{minipage}		
\begin{minipage}[b]{0.23\linewidth}
	\centering
	 \includegraphics[width=3.2cm]{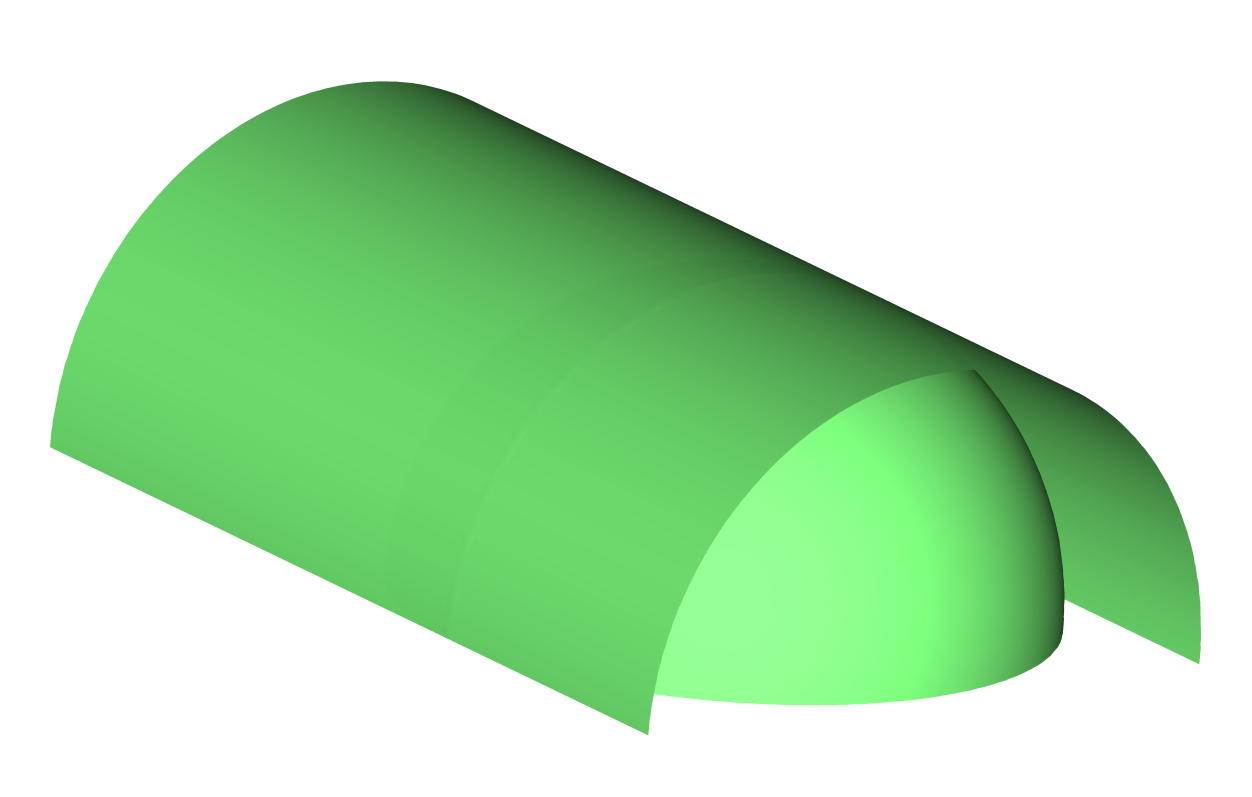}{\begin{center}(d) quartic approximation\end{center}}\label{fig:quartic}
\end{minipage}
\end{tabular}
\caption[Simultaneous approximation of triangular and tensor-product patches]{Implicit approximations of the surface described in Section \ref{subsec:several}. 
Note that the quartic approximation, which is an exact implicitization up to rounding error, is defined by the product of two polynomials and hence extra branches are present.}\label{fig:patches}						
\end{figure*}
In many circumstances it may be desirable to approximate several surface patches simultaneously, by a single implicit surface. This is possible using either the original or weak methods \cite{dokken_thesis,dokken_weak}. 

Suppose we have several parametric surfaces $\mathbf{p}_1(\mathbf{s}),\ldots,\mathbf{p}_r(\mathbf{s}).$ To find an implicit surface that approximates all these surfaces we may proceed as before to build matrices $\mathbf{D}_i$ corresponding to the individual manifolds $\mathbf{p}_i(\mathbf{s}).$ However, before performing the SVD, we stack the matrices to define \[\mathbf{D} = \left(\begin{smallmatrix} \mathbf{D}_1 \\ \vdots \\ \mathbf{D}_r \end{smallmatrix}\right).\] 

Using the weak form we build matrices $\mathbf{M}_i$ corresponding to the manifolds $\mathbf{p}_i(\mathbf{s}),$ but instead of stacking, we sum the matrices to form \[\mathbf{M} = \sum_{i=1}^r \mathbf{M}_i.\] Performing an SVD on $\mathbf{M}$ then defines the weak approximation. 

In fact, we are not restricted to approximating surfaces of the same type. There is also the possibility to simultaneously approximate points, curves and surfaces with different parametric forms. To exemplify this we approximate a surface defined by two rational tensor-product B\'ezier patches describing a half-cylinder, and two rational B\'ezier triangles describing a quarter-sphere, as pictured in Figure \ref{fig:patches}(a). The quadratic, cubic and quartic approximations displayed in Figure \ref{fig:patches} demonstrate some interesting properties of implicit representations. The quadratic approximation, which is in fact described by an ellipsoid, is clearly quite different from the parametric surface and, for most purposes, would not be a sufficient approximation. The cubic approximation is visually what we expect to see, and does indeed provide a close approximation. When we increase the degree to four, as expected, we obtain an exact implicit representation, though this is defined by the product of two polynomials which describe the cylindrical and spherical parts separately. Consequently, when visualizing the surface we see extra branches that are not present in the parametric representation.

\subsection{Approximate Implicitization of Surfaces with Singularities}\label{sec:bez2}

A simple example of a quadratic B\'ezier triangle $\mathbf{p}_2(\mathbf{s})$ that exhibits singularities is constructed by taking the three corner control points to be at the Cartesian origin $(0,0,0)$, and the three central control points to be at (1,0,0), (0,1,0) and (0,0,1). This is pictured in Figure \ref{fig:beziertri}(b). An exact quartic implicitization and an approximate cubic implicitization of $\mathbf{p}_2(\mathbf{s})$ are pictured in Figure \ref{fig:first}. We will now compare the approximations of this example with the example from Section \ref{subsec:firstex}, to see how the singular surface suffers from worse approximations. Table 1 lists the singular values for implicit approximations up to degree four, obtained by the original method. Both of the surfaces require degree four for an exact implicitization. However, the singular values of $\mathbf{p}_1(\mathbf{s})$ are much smaller for the quadratic and cubic approximations, indicating better approximations.
\begin{table}[t]
\begin{center}
 \begin{tabular}{ | c | c | c | c | c | c | }
  \hline
  Degree $m$ & 1 & 2 & 3 & 4 \\ \hline 
  $\sigma_{\min}$ of $\mathbf{p}_1(\mathbf{s})$ & 1.0 & 0.22984 & 0.047868 & 0.0 \\ \hline
  $\sigma_{\min}$ of $\mathbf{p}_2(\mathbf{s})$ & 1.0 & 0.62773 & 0.31596 & 0.0 \\ 
    \hline
 \end{tabular}\label{tab:comp}
\end{center}
\medskip { Table 1.} {Difference in the smallest singular values of $\mathbf{D}$ for a B\'ezier triangle with singularities $(\mathbf{p}_2(\mathbf{s}))$ and without singularities $(\mathbf{p}_1(\mathbf{s})).$}
\end{table}

\section{Conclusion}

This paper described how the original and weak methods of approximate implicitization can be applied to triangular B\'ezier surfaces. It presented examples which exhibit properties of the various approaches to approximate implicitization. It also highlighted ways in which to improve the efficiency of the algorithm in the numerical case, by exploiting symmetries in the calculations.

\section*{Acknowledgements}

This work has been supported by the European Community under the Marie Curie Initial Training Network ``SAGA - Shapes, Geometry and Algebra''  Grant Agreement Number 21458, and by the Research Council of Norway through the IS-TOPP program. We would like to thank Johan Simon Seland at SINTEF for making his real-time algebraic surface visualization system available to us (Figure \ref{fig:first}). We have also made use of the software Axel, \url{http://axel.inria.fr/} (Figures \ref{fig:beziertri} and \ref{fig:patches}).

\addcontentsline{toc}{section}{Bibliography}
\nocite{*}
\bibliographystyle{plain}
\bibliography{aiotbs}

\end{document}